\documentclass[a4paper,10pt]{article}

\usepackage[utf8x]{inputenc}

\title{Bounds On Factors Of Odd Perfect Numbers}
\author{Siddhartha Basak \\ St. Xaviers' Collegiate School, Kolkata - 700016, India \\ siddhartha.basak2012@gmail.com}

\usepackage{epsfig}
\usepackage{graphics}

\usepackage{amsmath,amssymb, amsthm}
\usepackage{color}

\newcommand{\bc}{\begin{center}}
\newcommand{\ec}{\end{center}}
\newcommand{\bq}{\begin{quote}}
\newcommand{\eq}{\end{quote}}
\newcommand{\btab}{\begin{tabular}}
\newcommand{\etab}{\end{tabular}}

\newcommand{\be}{\begin{equation}}
\newcommand{\ee}{\end{equation}}
\newcommand{\beqa}{\begin{eqnarray*}}
\newcommand{\eeqa}{\end{eqnarray*}}
\newcommand{\beqn}{\begin{eqnarray}}
\newcommand{\eeqn}{\end{eqnarray}}
\newcommand{\bbibl}{\begin{thebibliography}{99}}
\newcommand{\ebibl}{\end{thebibliography}}

\newcommand{\nn}{\nonumber}

\newcommand{\ba}{\begin{array}}
\newcommand{\ea}{\end{array}}

\renewcommand{\qed}{\rule{1.3mm}{3mm}}

\newcounter{cnt1}
\newcounter{cnt2}
\newcounter{cnt3}
\newcommand{\blr}{\begin{list}{$($\roman{cnt1}$)$} {\usecounter{cnt1}
                \setlength{\topsep}{0pt} \setlength{\itemsep}{0pt}}}
\newcommand{\bla}{\begin{list}{$($\alph{cnt2}$)$} {\usecounter{cnt2}
                \setlength{\topsep}{0pt} \setlength{\itemsep}{0pt}}}
\newcommand{\bln}{\begin{list}{$($\arabic{cnt3}$)$} {\usecounter{cnt3}
                \setlength{\topsep}{0pt} \setlength{\itemsep}{0pt}}}
\newcommand{\el}{\end{list}}

\newtheorem{thm}{\indent{\sc Theorem}}[section]

\newtheorem{rem}{\indent{\bf Remark}}[section]

\date{November 17, 2012}

\begin{document}
\maketitle
\begin{abstract}
 Much recent progress has been made concerning the probable existence of Odd Perfect Numbers, forming part of what has come to be known
 as Sylvester's Web Of Conditions \cite{1}. This paper proves some results concerning certain properties of the sums of reciprocals of the factors
 of odd perfect numbers, or, in more technical terms, the properties of the subsums of $\sigma_{-1} (n)$. By this result, it also 
 establishes strong bounds on the prime factors of odd perfect numbers using the number of distinct prime factors it may possess.
\end{abstract}

\section{Introduction}
Euclid himself, in book IX of his magnum opus, The Elements, stated and proved a method for finding even perfect numbers. Many hundred years later,
Euler proved that this method found all even perfect numbers, though his method said nothing about the odd kind. The earliest references to the 
Odd Perfect Number Quandary (the problem of proving the existence/non-existence of odd perfect numbers) can be found in the mathematical letters 
between Fr. Marin Mersenne and Rene Descartes in 1638, in which Descartes proposed that these elusive entities might, indeed, exist. 
Generations of mathematicians, amateur or professional, have attacked this problem \cite{1}. \\
{\bf The Eulerian Form:}  The Eulerian Form \cite{10} of odd perfect numbers (referred to as EF throughout the paper) comprises, arguably, the most 
important contribution to Sylvester's Web Of Conditions, due to Euler. It states that if a number n is an 
odd perfect number (OPN), then n is of the form $n = {p^b} {q_1}^{2a_1} {q_2}^{2a_2} \dots {q_r}^{2a_r}$ where $p, q_1, q_2, \dots q_r$ are prime, and 
$p\equiv b \equiv 1 \ \ (mod \ \ 4)$. \\

\medskip
The paper is organised as follows: \\ 
As a preliminary, we prove a well-known result concerning $\sigma_{-1}(n)$ (the sum of the reciprocals of the factors of n). \\
We establish, in Theorem {\ref{mt}}, stringent bounds on the reciprocal factor sum of an OPN, irrespective of the exponents
of the primes in its prime factorisation. \\
In Theorem {\ref{mt2}}, we establish an upper bound on the prime factors of n, using the previous theorem and the number of distinct prime factors
of the OPN. We provide tables of these upper bounds.

\section{A Preliminary}
All through the paper, n is an OPN and p a prime.
We first, for the sake of completeness, give proof of a well-known result: \\
\begin{equation}
\label{p}
\sigma_{-1}(n) = \prod_{p | n} \sum_{k=0}^{h_p}  p^{-k} = 2
\end{equation}
where $h_p$ is the degree of p in the prime factorisation of n.

{\bf Proof of (\ref{p}).} 
We know, from the definition of odd perfect numbers, that $\sigma(n) = \prod_{p | n} \sum_{k=0}^{h_p}  p^k = 2n$, where $\sigma(n)$ is the divisor sum. 
Hence 
\begin{eqnarray}
2 = 
\frac{\sigma(n)}{n} 
&=& \frac{\prod_{p | n} \sum_{k=0}^{h_p}  p^k }{n} \nn\\
&=& \frac{\prod_{p | n} \sum_{k=0}^{h_p}  p^k }{\prod_{p | n}  p^{h_p} } \nn\\
&=& \prod_{p | n} \frac{ \sum_{k=0}^{h_p}  p^k }{p^{h_p} } \nn\\
&=& \prod_{p | n} \sum_{k=0}^{h_p}  p^{-k} 
\end{eqnarray}
\vspace{-10pt}
\hfill  \qed \\

We use this result in the main section of this paper.

\section{Main results}
 We start by proving bounds for
 $\prod_{p | n} \sum_{i=0}^{\alpha} p^{-i}$ for a predetermined $\alpha$.
\begin{thm}
\label{mt}
$$\frac{2^{\alpha + 2}}{\zeta(\alpha + 1) (2^{\alpha + 1}-1)} \ \ < \ \ \prod_{p | n} \sum_{i=0}^{\alpha} p^{-i} \ \ < \ \ 2$$
for a positive integer $\alpha \le h_p \ \forall p$ (where $\zeta(s)$ is the Riemann Zeta Function).
\end{thm}

{\bf Proof.}
We have observed above,
\begin{equation}
 \prod_{p | n} \sum_{k=0}^{h_p} p^{-k} = 2. \nn\\
\end{equation}
\vspace{-10pt}
 Let $h_p \ge \alpha \ \ \forall p$, for some $\alpha$.
 Then, 
\begin{eqnarray}
\label{lb}
2 \ = \prod_{p | n} \sum_{k=0}^{h_p} p^{-k} & \le & \prod_{p | n} \sum_{k=0}^{(\alpha + 1)\lfloor\frac{h_p}{\alpha + 1}\rfloor + \alpha} p^{-k} \ \ \ \ \ \ \ \ \ \ (\mbox{since $(\alpha + 1) \lfloor \frac{h_p}{\alpha + 1} \rfloor + \alpha \ge h_p$}) \nn\\
&=&  \prod_{p | n} \ \sum_{i = 0}^{\alpha} \sum_{j = 0}^{\lfloor \frac{h_p}{\alpha + 1} \rfloor} p^{-(j(\alpha + 1) + i)} \nn\\
&=&  \prod_{p | n} \Bigl({\sum_{i=0}^{\alpha} p^{-i} \sum_{j=0}^{\lfloor\frac{h_p}{\alpha + 1}\rfloor} p^{-j(\alpha + 1)}}\Bigr) \nn\\
&=&  \prod_{p | n} \sum_{i=0}^{\alpha} p^{-i} \prod_{p | n} \sum_{j=0}^{\lfloor\frac{h_p}{\alpha + 1}\rfloor} p^{-j(\alpha + 1)} \nn\\
&<&  \prod_{p | n} \sum_{i=0}^{\alpha} p^{-i} \prod_{p | n} \sum_{j=0}^{\infty} p^{-j(\alpha + 1)} \nn\\
&<&  \prod_{p | n} \sum_{i=0}^{\alpha} p^{-i} \prod_{p \ prime, \ p\neq2} \sum_{j=0}^{\infty} p^{-j(\alpha + 1)} \nn\\
&=&  \prod_{p | n} \sum_{i=0}^{\alpha} p^{-i} \frac{\prod_{p \ prime} \sum_{j=0}^{\infty} p^{-j(\alpha + 1)}}{\sum_{j=0}^{\infty} 2^{-j(\alpha + 1)}} \nn\\
&=&  \prod_{p | n} \sum_{i=0}^{\alpha} p^{-i} \frac{\prod_{p \ prime} \frac{1}{1-p^{-(\alpha + 1)}}}{\frac{1}{1-2^{-(\alpha + 1)}}} \nn\\
&=&  \Bigl(\zeta (\alpha + 1)/ \frac{2^{\alpha + 1}}{2^{\alpha + 1}-1}\Bigr) \prod_{p | n} \sum_{i=0}^{\alpha} p^{-i}  \ \ \ \ (since \ \zeta(s)=\prod_{p \ prime} \frac{1}{1-p^{-s}}) \nn\\
& \Rightarrow &  2  <  \Bigl(\zeta (\alpha + 1) \frac{2^{\alpha + 1}-1}{2^{\alpha + 1}}\Bigr) \prod_{p | n} \sum_{i=0}^{\alpha} p^{-i} \nn\\ 
& \Rightarrow &  \frac{2^{\alpha + 2}}{\zeta(\alpha + 1) (2^{\alpha + 1}-1)} < \prod_{p | n} \sum_{i=0}^{\alpha} p^{-i}.
\end{eqnarray}
And $\alpha \le h_p$. But we know, from the EF, that every exponent (the $h_p$s), is even other than the first. So the final product in (\ref{lb}) 
does not have $\alpha = h_p$ for all p.
\begin{equation}
 \therefore \ \ \prod_{p | n} \sum_{i=0}^{\alpha} p^{-i} < \prod_{p | n} \sum_{i=0}^{h_p} p^{-i} = 2 \nn\\
\end{equation}
To summarise,
\begin{equation}
\frac{2^{\alpha + 2}}{\zeta(\alpha + 1) (2^{\alpha + 1}-1)} \ \ < \ \ \prod_{p | n} \sum_{i=0}^{\alpha} p^{-i} \ \ < \ \ 2
\end{equation}
\hfill \qed \\

Also, since $\alpha = 1$ works ($h_p$ obviously $\geq 1$), we have
\begin{equation}
 \frac{8}{3\zeta(2)} = \frac{16}{\pi^2} < \prod_{p | n} \Bigl(1 + \frac{1}{p} \Bigr) < 2
\end{equation}

\begin{rem}
\label{cor}
 Note, however, that $\frac{16}{\pi^2} \approx 1.621138938$ is not as strong as the bound obtained by placing $\alpha = 2$, i.e., $\frac{16}{7\zeta(3)} 
\approx 1.901502566$ (computed on Wolfram Alpha \cite{2}). The $\alpha = 2$ case is plausible, since all other exponents besides the first are even, so
$\ge 2$, and the first is of the form 4k + 1 (see Introduction, Eulerian Form). This gives us two cases, (1) the first exponent, b, is 1, or (2) $b > 2$ 
(5, 9, ...). \\
Following the proof of Theorem ${\ref{mt}}$, and isolating the first prime, q, in the first line, we have \\
{\bf Case 1.}
\begin{equation}
\frac{16q^3}{7\zeta(3)(q^3-1)} < \Bigl(1 + \frac{1}{q} \Bigr) \prod_{p | n, p \neq q} \Bigl(1 + \frac{1}{p} + \frac{1}{p^2}\Bigr) < 2 \nn\\
\end{equation}
But since \\
\begin{equation}
\frac{16}{7\zeta(3)} < \frac{16q^3}{7\zeta(3)(q^3-1)}, \nn\\
\end{equation}

\begin{equation}
\frac{16}{7\zeta(3)} < \Bigl(1 + \frac{1}{q} \Bigr) \prod_{p | n, p \neq q} \Bigl(1 + \frac{1}{p} + \frac{1}{p^2}\Bigr) < 2
\end{equation}
The upperbound, 2, holds since the product still cannot have $\alpha = h_p \ \forall p$, since all the even exponents cannot be 2 \cite{6}. And the other 
case is simply $\alpha = 2$, i.e., \\
{\bf Case 2.}
\begin{equation}
\frac{16}{7\zeta(3)} < \prod_{p | n} \Bigl(1 + \frac{1}{p} + \frac{1}{p^2}\Bigr) < 2
\end{equation}
This provides an efficient, though not optimal method for proving a number is NOT an OPN, since one does not need to know the exponents of the primes (thus
not needing a complete factorisation). These two cases also make evident why all three of 3,5 and 7 cannot divide an OPN.
\end{rem}
\bigskip

\begin{thm}
\label{mt2}
If $n$ is an OPN and $\omega(n) = m$ is the number of distinct prime factors of n, then there exist primes $p_{I_1}, \ p_{I_2} \ and \ p_{I_3}$ such that
the first, second and third prime factors are less than the respective $p_{I_k}$s, where the $p_{I_k}$s can be determined, given $m$, 
($p_r$ is the $r^{th}$ prime).
\end{thm}

{\bf Proof.}
Let $n$ (arbitrary odd perfect number) be written as,
\begin{equation}
\label{eq21}
 n = p_{i_1}^{a_1}  p_{i_2}^{a_2} \cdots  p_{i_m}^{a_m}. \nn\\
\end{equation}

We rewrite the product from Theorem ({\ref{mt}}) as
\begin{equation}
\label{eq22}
 \prod_{p_{i_j} | n, \ j=1 }^{j=m} \Bigl(1 + \frac{1}{p_{i_j}} \Bigr) 
\end{equation}
and define
\begin{equation}
\label{eq23}
\rho^{(1)}_r= \prod_{j=r}^{r+m-1} \Bigl(1 + \frac{1}{p_{j}} \Bigr) . 
\end{equation}
It is evident that $\rho^{(1)}_{i_1} \ge$ the product of ({\ref{eq22}}). But note that
\begin{equation}
\label{eq24}
\lim_{r \rightarrow \infty} \rho^{(1)}_r = 1 < \frac{16}{\pi^2}  . \nn\\
\end{equation}
Therefore, \\
\begin{eqnarray}
\label{eq25}
\exists \ I_1 \ such \ that \ \rho^{(1)}_r < \frac{16}{\pi^2} \ \ \forall r > I_1 \nn\\ 
\therefore \ \ i_1 \le I_1 \ and \ p_{i_1} \le p_{I_1}, \ if \ n \ is \ an \ OPN, \ as \ claimed.
\end{eqnarray}

Similarly, we define
\begin{eqnarray}
\label{eq26}
&& \rho^{(2)}_r= \Bigl(1+\frac{1}{3} \Bigr) \prod_{j=r}^{r+m-2} \Bigl(1 + \frac{1}{p_{j}} \Bigr), \ so \ that \nn\\
&& \rho^{(2)}_{i_2} \ge \Bigl(1+ \frac{1}{p_{i_1}} \Bigr) \prod_{p_{i_j} | n, \ j=2 }^{j=m} \Bigl(1 + \frac{1}{p_{i_j}} \Bigr)
\end{eqnarray}
and note that
\begin{eqnarray}
\label{eq27}
\lim_{r \rightarrow \infty} \rho^{(2)}_r = 1+\frac{1}{3} < \frac{16}{\pi^2}  . \nn\\
\therefore \ \exists \ I_2 \ such \ that \ \rho^{(2)}_r < \frac{16}{\pi^2} \ \  \forall r > I_2 \nn\\ 
\therefore \ \ i_2 \le I_2 \ and \ p_{i_2} \le p_{I_2}, \ if \ n \ is \ an \ OPN, \ as \ claimed.
\end{eqnarray}

And finally,
\begin{eqnarray}
\label{eq28}
&& \rho^{(3)}_r= \Bigl(1+\frac{1}{3} \Bigr) \Bigl(1+\frac{1}{5} \Bigr) \prod_{j=r}^{r+m-3} \Bigl(1 + \frac{1}{p_{j}} \Bigr), \ so \ that \nn\\
&& \rho^{(3)}_{i_3} \ge \Bigl(1+ \frac{1}{p_{i_1}} \Bigr) \Bigl(1+ \frac{1}{p_{i_2}} \Bigr) \prod_{p_{i_j} | n, \ j=3 }^{j=m} \Bigl(1 + \frac{1}{p_{i_j}} \Bigr)
\end{eqnarray}

\begin{eqnarray}
\label{eq29}
\lim_{r \rightarrow \infty} \rho^{(3)}_r = \Bigl(1+\frac{1}{3} \Bigr) \Bigl(1+\frac{1}{5} \Bigr) < \frac{16}{\pi^2}  . \nn\\
\therefore \ \exists \ I_3 \ such \ that \ \rho^{(3)}_r < \frac{16}{\pi^2} \ \  \forall r > I_3 \nn\\ 
\therefore \ \ i_3 \le I_3 \ and \ p_{i_3} \le p_{I_3}, \ if \ n \ is \ an \ OPN, \ as \ claimed.
\end{eqnarray}

\hfill  \qed \\

Unfortunately, $\Bigl(1+\frac{1}{3} \Bigr) \Bigl(1+\frac{1}{5} \Bigr) \Bigl(1+\frac{1}{7} \Bigr) > \frac{16}{\pi^2}$, so these methods do not apply to
other prime factors.

\medskip 
A table of the values of $p_{I_k}$ for each $m$ upto 20, follows. The list begins from 9, since $\omega(n) \ge 9$ if $n$ is an OPN, as established by Nielsen
in \cite{4}.\\

\medskip
\begin{center}
\begin{tabular}{|c|c|c|c|}
\hline
& \multicolumn{2}{c}{$p_{I_k}$ (for $\alpha=1$)} & \\ \hline
m & k=1 & k=2 & k=3 \\ \hline
9 & 11 & 31 & 509 \\ \hline
10 & 11 & 31 & 593 \\ \hline
11 & 11 & 37 & 659 \\ \hline
12 & 13 & 41 & 739 \\ \hline
13 & 13 & 43 & 811 \\ \hline
14 & 13 & 43 & 881 \\ \hline
15 & 13 & 47 & 947 \\ \hline
16 & 13 & 53 & 1031 \\ \hline
17 & 17 & 53 & 1093 \\ \hline
18 & 17 & 59 & 1171 \\ \hline
19 & 17 & 61 & 1237 \\ \hline
20 & 17 & 61 & 1301 \\ \hline
\end{tabular}
\end{center}

\medskip 
 A well-known result due to Perisastri \cite{5} states that an OPN with $m$ distinct prime factors has its lowest prime factor $\le \frac{2}{3}m + 3$. The 
method proved in this paper produces marginally stronger bounds.

\medskip
These results, combined with previous results concerning the prime factors of OPNs, may be used for future research into the problem.

\section{Acknowledgments}
I thank my father, G. K. Basak, for many insightful discussions on this topic, and for much integral assistance with this paper.

\end{document}